\documentclass[danish, a4paper, 12pt]{report}
\usepackage[paperheight=28cm, paperwidth=20cm, margin=2cm]{geometry}
\usepackage[T1]{fontenc}
\usepackage[utopia]{mathdesign}
\usepackage[english]{babel}
\usepackage[utf8]{inputenc}
\usepackage[usenames,dvipsnames]{color}
\usepackage{amsmath}
\usepackage{graphicx}
\usepackage[toc,page]{appendix}
\usepackage{multicol}
\usepackage{tabularx}
\usepackage{amsopn}
\usepackage{xparse}
\usepackage{booktabs}
\usepackage{pgfplotstable}
\usepackage{rotating}
\usepackage{lscape}
\usepackage{wrapfig}
\usepackage{verbatim}
\usepackage{physics}
\usepackage{wasysym}
\usepackage{pstricks}
\usepackage{pspicture}
\usepackage{mathrsfs}
\graphicspath{{"Immagini/"}}
\usepackage{mathtools}

\DeclarePairedDelimiter\floor{\lfloor}{\rfloor}
\definecolor{ballblue}{rgb}{0.13, 0.67, 0.8}
\definecolor{apple}{rgb}{0.55, 0.71, 0.0}
\definecolor{awesome}{rgb}{1.0, 0.13, 0.32}
\definecolor{azure}{rgb}{0.0, 0.5, 1.0}
\DeclareMathOperator{\de}{\text{d}}
{\left\lbrace\begin{array}{@{}l@{}}}%
{\end{array}\right.}

\author{ Lectures from a seminar of Mathematical Physics by Enrico Masina$\dagger$  \\ \\ \textsf{$\dagger$ Dipartimento di Fisica ed Astronomia (DIFA), University of Bologna} \\ \textsf{ Alma Mater Studiorum,} \textsf{ and INFN, Italy} \\ \texttt{\textcolor{blue}{Enrico.Masina@bo.infn.it}} \\\\}
\date{}
\title{\textbf{\textrm{\textbf{A review on the Exponential-Integral special function and other strictly related special functions.\\}}}}

\begin{document}
\maketitle 
\noindent\textbf{\huge Introduction}
\\\\
The aim of this short series of lectures is to provide a good treatise about the Exponential Integral special function, covering also a few selected topics like its generalisation, other related special functions (Sine Integral, Cosine Integral and the Logarithmic Integral) and the theory of Asymptotic Series.\\
Special functions have always played a central role in physics and in mathematics, arising as solutions of particular differential equations, or integrals, during the study of particular important physical models and theories in Quantum Mechanics, Thermodynamics, Fluid Dynamics, Viscoelasticity, et cetera. \\
The theory of Special Functions is closely related to the theory of Lie groups and Lie algebras, as well as certain topics in Mathematical Physics, hence their use and a deep knowledge about them is as necessary as mandatory.\\\\
In its first part, the present paper aims to give a particular detailed treatise over the Exponential Integral function and strictly related functions: its generalisation (the so called generalised-Ei), and the Modified Exponential Integral (often named Ein), giving their Series expansion and showing their plots over the real plane.\\\\
The second part is entirely dedicated to a review of Asymptotic Series, starting with the Big-O and Little-o symbols (the Landau symbols), in order to reconnect to the Exponential Integral (and related functions) to study their asymptotic behaviour. The end of this part will also look at what happens when the argument of the Exp-Integral becomes imaginary whence the study of the Sine Integral (Si) and Cosine Integral (Ci) special functions.\\\\
The third part will focus on another special function strictly connected to Ei: the Logarithmic Integral providing also a brief recall of its connection with the Prime Numbers counting function $\pi(x)$.\\\\
The last part will naively present a solved physical problem where the Exponential Integral pops out, eventually leaving to the readers three good exercises in order to have some practice and become familiar with those important special functions.

\tableofcontents

\chapter{The Exponential Integral Function}

We can immediately start by giving the mathematical expression of the so called \textbf{Exponential Integral Function}, which is defined as
$$\boxed{\text{Ei}(x) = \int_{-\infty}^x \frac{e^t}{t}\ \text{d}t} ~~~~~~~~~~~ x\in\mathbb{W} \footnote{The set $\mathbb{W}$ is defined as the \textbf{Whole} set, namely the set of all the natural numbers but zero:
$$\mathbb{W} \equiv \mathbb{R}^+\backslash\{0\} = \{1, 2, 3, 4, \ldots \}$$
We use instead the notation $\mathbb{R}^+$ to indicate the set of the natural numbers plus zero:
$$\mathbb{R}^+ = \{0, 1, 2, 3, \ldots \}$$}$$
So we can see the Exponential Integral is defined as one particular definite integral of the ratio between the exponential function and its argument.\\\\
The Risch Algorithm can be performed to show that this integral is not an elementary function, namely a primitive of $\text{Ei}(x)$ in terms of elementary functions does not exist.\\
Such a function has a pole at $t = 0$, hence we shall interpret this integral as the \textit{Cauchy Principal Value}:\\
$$\text{Ei}(x) = \lim_{\alpha\to 0} \left[\int_{-\infty}^{\alpha} \frac{e^t}{t}\ \text{d}t + \int_{\alpha}^x \frac{e^t}{t}\ \text{d}t\right]$$
\\
We can define $\text{Ei}(x)$ in a more suitable way through a parity transformation
$$\begin{cases}
t\to -t \\
x\to -x
\end{cases}
$$
one gets
$$\boxed{\mathcal{E}_1(x) = -\text{Ei}(-x) = \int_x^{+\infty} \frac{e^{-t}}{t}\ \text{d}t}$$
\\\\
When the argument takes complex values, the definition of the integral becomes ambiguous due to branch points at $0$ and $\infty$. \\
Introducing the complex variable $z = x + iy$, we can use the following notation to define the \textbf{Exponential Integral in the Complex plane}\footnote{A brief summary about complex numbers: \\
Given a complex number $z = x+iy$ we can always write it in the exponential form
$$z = |z| e^{i\theta}$$
where
$$|z| = \sqrt{x^2 + y^2} ~~~~~~~~~~~ \theta = \arg(z) = \arctan\left(\frac{y}{x}\right)$$}:
$$\boxed{\mathcal{E}_1(z) = \int_z^{+\infty} \frac{e^{-t}}{t}\ \text{d}t} ~~~~~~~~~~~ |\arg(z) < \pi |$$
\\
To clarify the graphical behaviour of the two functions, the following plot may come in handy.
\\
\begin{figure}[h!]
\centering
\includegraphics[scale=1.7]{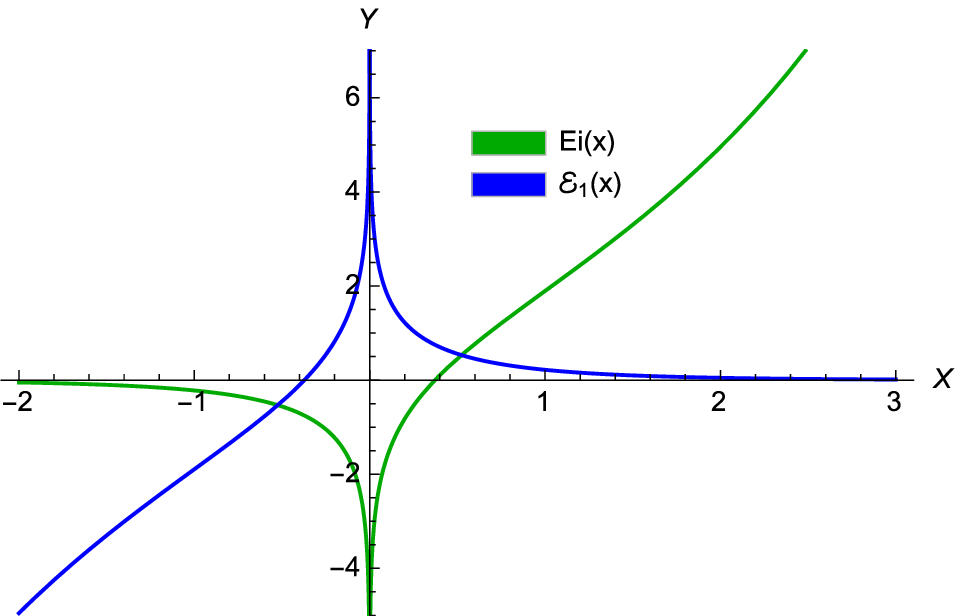}
\end{figure}
\
\\
We can immediately write down some useful known values:
$$\text{Ei}(0) = -\infty ~~~~~~~~~~~~~ \text{Ei}(-\infty) = 0 ~~~~~~~~~~~~~ \text{Ei}(+\infty) = +\infty$$
$$\mathcal{E}_1(0) = +\infty ~~~~~~~~~~~~~ \mathcal{E}_1(+\infty) = 0 ~~~~~~~~~~~~~ \mathcal{E}_1(-\infty) = -\infty$$
It's actually simple to find the values of $\mathcal{E}_1(x)$ from $\text{Ei}(x)$ (and vice versa) by using the previously written relation: 
$$\mathcal{E}_1(x) = -\text{Ei}(-x)$$
\\
We notice that the function $\mathcal{E}_1(x)$ is a monotonically decreasing function in the range $(0; +\infty)$.\\
The function $\mathcal{E}_1(z)$ is actually nothing but the so called \textbf{Incomplete Gamma Function}: \footnote{Let's recall the definition of the \textbf{Ordinary Gamma Function}
$$\Gamma(s) = \Gamma(s, 0) = \int_0^{+\infty} t^{s-1} e^{-t}\ \text{d}t$$
whilst
$$\Gamma(0, 0) = \int_0^{+\infty} \frac{e^{-t}}{t}\ \text{d}t = \infty$$
}
$$\mathcal{E}_1(z) \equiv \Gamma(0, z)$$
where
$$\Gamma(s, z) = \int_z^{+\infty} t^{s-1} e^{-t}\ \text{d}t$$
Indeed by putting $s = 0$ we immediately find $\mathcal{E}_1(z)$.
\\\\
By introducing the \textbf{small Incomplete Gamma Function}
$$\gamma(s, z) = \int_0^z t^{s-1} e^{-t}\ \text{d}t$$
we are able to write down a very obvious, sometimes useful and straightforward relation between the three Gamma functions:
$$\Gamma(s, 0) = \gamma(s, z) + \Gamma(s, z)$$
\\\\
Let's now come back to the Exponential Integral. \\
Let's perform a naive change of variable
$$t \to zu ~~~~~~~ \text{d}t = z\ \text{d}u$$
Step by step we get:
$$\int_z^{+\infty} \frac{e^{-t}}{t}\ \text{d}t \to \int_1^{+\infty} \frac{e^{-zu}}{zu}z\ \text{d}u = \int_1^{+\infty} \frac{e^{-zu}}{u}\ \text{d}u$$
In this way we define the \textbf{General Exponential Integral}:
$$\boxed{\mathcal{E}_n(z) = \int_1^{+\infty} \frac{e^{-zu}}{u^n}\ \text{d}u ~~~~~~~~~~~ n\in\mathbb{R}}$$
with the particular value
$$\mathcal{E}_n(0) = \frac{1}{n-1}$$
\newpage
\noindent
We can show the behaviour of the first five functions $\mathcal{E}_n(x)$, namely for $n = 0, \ldots 5$:
\\
\begin{figure}[h!]
\centering
\includegraphics[scale=1.5]{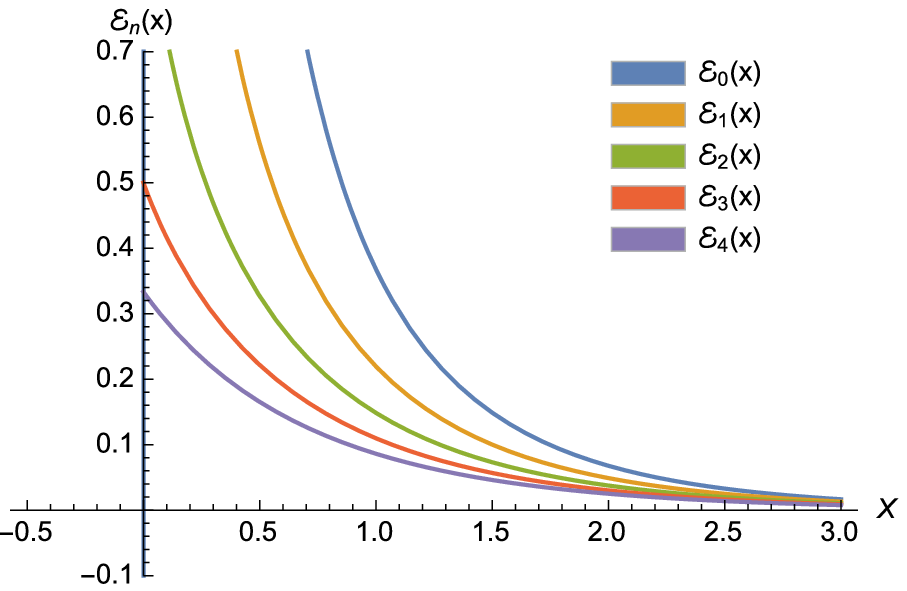}
\end{figure}
\
\\
\section{Series Expansion for $\text{Ei}(z)$}

Let's go back to the Exponential Integral we defined at the very beginning of the paper.
$$\text{Ei}(z) = \int_{-\infty}^z \frac{e^t}{t}\ \text{d}t$$
What we are going to do is to split the integral initially into three pieces:
$$\int_{-\infty}^z \frac{e^t}{t}\ \text{d}t = \int_{-\infty}^{-1} \frac{e^t}{t}\ \text{d}t + \textcolor{Red}{\int_{-1}^0 \frac{e^t}{t}\ \text{d}t} + \textcolor{cyan}{\int_{0}^z \frac{e^t}{t}\ \text{d}t}$$
and now with a mathematical trick we add and cut the term
$$\int_{-1}^z \frac{\text{d}t}{t}$$
as follows:
$$\int_{-\infty}^{-1}\frac{e^t}{t}\ \text{d}t + \int_{-1}^0 \frac{e^t}{t}\ \text{d}t + \int_{0}^z \frac{e^t}{t}\ \text{d}t \textcolor{apple}{+ \int_{-1}^z \frac{\text{d}t}{t} - \int_{-1}^z \frac{\text{d}t}{t}}$$
Expanding the last integral into $\int_{-1}^0 + \int_0^z$ we get
$$\int_{-\infty}^{-1}\frac{e^t}{t}\ \text{d}t + \int_{-1}^0 \frac{e^t}{t}\ \text{d}t + \int_{0}^z \frac{e^t}{t}\ \text{d}t \textcolor{apple}{+ \int_{-1}^z \frac{\text{d}t}{t} - \int_{-1}^0 \frac{\text{d}t}{t} - \int_{0}^z \frac{\text{d}t}{t}}$$
This appears as a weird mathematical manipulation but now we can add up the similar integral, namely the second term from the left with the second to last term, and the third term from the left with the very last term at the right side, getting:
$$\int_{-\infty}^{-1}\frac{e^t}{t}\ \text{d}t +\int_{-1}^0 \frac{e^t-1}{t}\ \text{d}t + \int_{0}^z \frac{e^t-1}{t}\ \text{d}t + \int_{-1}^z \frac{\text{d}t}{t}$$
\\
Now it's about a bit of maths.
\\\\
Let's take the very first term, and let's perform the reciprocal transformation of variable $t \to 1/t$:
$$\int_{-\infty}^{-1} \frac{e^{t}}{t}\ \text{d}t \xrightarrow{~~ t \to 1/t ~~} \int_{-1}^0 \frac{e^{1/t}}{t}\ \text{d}t$$
We did this because now we can sum this term with the second term above which leads us to
$$\int_{-1}^0 \frac{e^{1/t} + e^t - 1}{t}\ \text{d}t$$
and again with a parity transformation $t \to -t$ we write it down in a well known form (this is a really important and interesting integral):
$$\int_0^1 \frac{1 - e^{-t} - e^{-1/t}}{t}\ \text{d}t = \gamma = -\psi(1)$$
where $\psi(z)$ is the logarithmic derivative of the Gamma function a.k.a. the \textbf{Digamma Function}:
$$\psi(z) = \frac{\text{d}}{\text{d}z}\ln\Gamma(z)$$
and $\gamma$ is the famous Euler-Mascheroni constant
$$\gamma \approx 0.5772156649015328606065\ldots $$
\hrule 
\subsection*{Brief digression over $\psi(z)$, $\gamma$ and $\Gamma(z)$}

Having introduced the Euler-Mascheroni constant, and the Digamma Function, it's necessary to show some connections between them.\\
First of all, let's recall the Gamma function and its derivative:
$$\Gamma(z) = \int_0^{+\infty} t^{z-1} e^{-t}\ \text{d}t$$
$$\Gamma'(z) = \frac{\text{d}}{\text{d}z} \Gamma(z) = \int_0^{+\infty} t^{z-1}e^{-t}\ln(t)\ \text{d}t$$
\\
We define the Digamma function as
$$\psi(z) = \frac{\Gamma'(z)}{\Gamma(z)} = \frac{\text{d}}{\text{d}z}\ln\Gamma(z)$$
\\
which is actually generally written as $\psi^0(z)$ because it belongs to a general class of functions called \textbf{Polygamma Functions}:
$$\psi^{m}(z) = \frac{\text{d}^{m+1}}{\text{d}z^{m+1}}\ln\Gamma(z)$$
in which we recognize the Digamma function $\psi^0(z)$ for $m = 0$.
\\\\
There is a well known value for the Digamma function, which is at $z = 1$ where
$$\psi(1) = -\gamma$$
Because of what we just wrote above, we have
$$\gamma = -\psi(1) = -\Gamma'(1) = -\int_0^{+\infty}e^{-t}\ln(t)\ \text{d}t$$
and by virtue of what has been found before, we have a first relationship between the Digamma function and the integral we found:
$$\int_0^1 \frac{1 - e^{-t} - e^{-1/t}}{t}\ \text{d}t = \gamma = -\int_0^{+\infty}e^{-t}\ln(t)\ \text{d}t$$
The integral in the left-hand side is not trivial, but the one on the right-hand side is, and we will give its solution in a while.
\\\\
If one took a general course in Series and Calculus, he won't be surprised in seeing the Euler-Mascheroni constant here since it already came out from a well known series, that is the Harmonic Series, and more precisely:
$$\gamma = \lim_{k\to +\infty} \left(\sum_{n = 1}^{k} \frac{1}{n} - \ln(k)\right)$$
That is easily numerically provable. Let's denote the bracket as $G(k)$. We find that
$$G(3) = \sum_{n = 1}^{3} \frac{1}{n} - \ln(3) = 0.734721\ldots $$
$$G(11) = \sum_{n = 1}^{11} \frac{1}{n} - \ln(11) = 0.621982\ldots $$
$$G(47) = \sum_{n = 1}^{47} \frac{1}{n} - \ln(47) = 0.587816\ldots $$
$$G(859) = \sum_{n = 1}^{859} \frac{1}{n} - \ln(859) = 0.577798\ldots $$
And going to infinity we find exactly the Euler-Mascheroni constant.
\\
\hrule
\newpage
\noindent
Now, let's go on with the Exponential Integral question. Aside the previous term we found, we are left with two integrals the very last one of which is trivial:
$$\int_{-1}^z \frac{\text{d}t}{t} = \ln(-z)$$
\\
Finally, adding all the terms together we get
$$\text{Ei}(z)  = \gamma + \ln(-z) + \int_0^z \frac{e^t - 1}{t}\ \text{d}t ~~~~~~~~~~~ |\arg(-z)| < \pi $$
The last integral can be now expanded in series since it does converge over the whole plane:
$$\int_0^z \frac{e^t - 1}{t}\ \text{d}t = \int_0^z \sum_{k = 1}^{+\infty} \frac{t^{k-1}}{k!}\ \text{d}t = \sum_{k = 1}^{+\infty} \frac{z^k}{k\cdot k!}$$
which finally leads us to the \textbf{\textcolor{awesome}{Series representation for the Exponential Integral}}:
$$\boxed{\text{Ei}(z) = \gamma + \ln(-z) + \sum_{k = 1}^{+\infty} \frac{z^k}{k\cdot k!} ~~~~~~~~~~~ |\arg(-z)| < \pi}$$
\\
From a very important relationship in the real field:
$$\mathcal{E}_1(-x \pm i0) = -\text{Ei}(x) \mp i\pi$$
we immediately get the \textbf{Series representation with a real argument}:
$$\text{Ei}(x) = \gamma  + \ln(x) + \sum_{k = 1}^{+\infty} \frac{x^k}{k\cdot k!}$$
\\
In the same way, recalling that $-\text{Ei}(-z) = \mathcal{E}_1(z)$ we can write down the \textbf{\textcolor{awesome}{Series representation for the Generalised Exponential Integral 1}}:
$$\boxed{\mathcal{E}_1(z) = -\gamma -\ln(z) - \sum_{k = 1}^{+\infty} \frac{(-z)^k}{k\cdot k!}}$$
which can also be written in this way, after a simple algebra:
$$\mathcal{E}_1(z) = -\gamma -\ln(z) + \sum_{k = 1}^{+\infty} \frac{(-1)^{k+1}(z)^k}{k\cdot k!}$$
\newpage\noindent 
\section*{Appendix A - Proof of: $-\int_0^{+\infty} e^{-t}\ln(t)\ \textsf{d} t   = \gamma$}

Recalling the limit definition of the exponential function
$$e^t = \lim_{k\to +\infty}\left(1 - \frac{t}{k}\right)^k  \de t$$
and substituting it into the integral, we get:
$$-\int_0^{+\infty} e^{-t}\ln(t)\ \textsf{d} t   =  -\int_0^{+\infty} \ln(t) \lim_{k\to +\infty}\left(1 - \frac{t}{k}\right)^k \de t = - \lim_{k\to +\infty}\int_0^{+\infty}  \left(1 - \frac{t}{k}\right)^k  \de t$$
Now we perform the change of variable 
$$u = 1 - \frac{t}{k} ~~~~~~~ t = k - ku ~~~~~~~ \de t = -k \de u$$
and the extrema of the integration vary from $1$ to $0$. \\
We have then
$$
\begin{array}{rll}
-\lim_{k\to +\infty} \int_0^1 ku^k \ln(k - ku)\ \de u & = -\lim_{k\to +\infty} k \int_0^1 u^k \left[\ln(k) - \ln(1-u)\right]\ \de u \\\\
& = -\lim_{k\to +\infty} k \left\{ \int_0^1 \left(u^k\ln(k) - u^k\ln(1-u)\right)\ \de u\right\}
\end{array}
$$
The two integration are quite easy to compute. A repeated integration by parts is required, and it will lead to a result in terms of a series.\\
We have then:
$$ -\lim_{k\to +\infty} k \left\{\frac{1}{k+1}\ln(k) - \sum_{n = 1}^k \frac{1}{n(n+k)}\right\}
$$
Now we make use of the well known identity (from Telescopic Series)
$$\sum_{n = 1}^k \frac{1}{n(n+k)} = \frac{1}{k}\sum_{n = 1}^k \frac{1}{n}$$
therefore
$$
\begin{array}{rll}
-\displaystyle\lim_{k\to +\infty} k \left\{\frac{1}{k+1}\ln(k) - \sum_{n = 1}^k \frac{1}{n(n+k)}\right\} & = -\displaystyle\lim_{k\to +\infty}\left(\frac{k}{k+1}\ln(k) - \frac{k}{k}\sum_{n = 1}^k \frac{1}{n}\right)\\\\
& = \displaystyle\lim_{k\to +\infty} -\frac{k\ln(k)}{k+1} + \sum_{n = 1}^k\frac{1}{n}\\\\
\text{and approaching $k$ to infinity} & = \displaystyle\lim_{k\to +\infty} -\ln(k) + \sum_{n = 1}^{k}\frac{1}{n} \\\\
& = \gamma
\end{array}
$$
As wanted.
\newpage\noindent 

\section*{Appendix B - Proof of: $\int_0^1 \frac{1 - e^{-t} - e^{-1/t}}{t}\ \textsf{d}t = \gamma$}

$$\int_0^1 \frac{1 - e^{-t} - e^{-1/t}}{t}\ \de t = \int_0^1 \frac{1 - e^{-t}}{t}\ \de t - \int_0^1 \frac{e^{-1/t}}{t}\ \de t$$
Now the first term comes integrated by parts, obtaining
$$\textcolor{awesome}{-\int_0^1 e^{t}\ln(t)\ \de t}$$
Whilst in the second term let's replace
$$t \to \frac{1}{t}$$
in order to get
$$-\left(- \int_{+\infty}^1 \frac{e^{-t}}{t}\ \de t\right) = \int_{+\infty}^1 \frac{e^{-t}}{t}\ \de t$$
which is now integrated by parts once to obtain
$$\textcolor{azure}{-\int_1^{+\infty}e^{-t}\ln(t)\ \de t}$$
Now we see the two arguments of the two integrals are identical, hence we can sum them in a single integral in $[0, +\infty)$:
$$-\int_0^1 e^{t}\ln(t)\ \de t - \int_1^{+\infty}e^{-t}\ln(t)\ \de t = - \int_0^{+\infty} e^{-t}\ln(t)\ \de t$$
In Appendix A we proved that the very last term in the equation above was exactly equal to $\gamma$, 
$$- \int_0^{+\infty} e^{-t}\ln(t)\ \de t = \gamma$$
As wanted.
\newpage\noindent 

\section{The Modified Exponential Integral :  $\text{Ein}(z)$}

The modified Exponential Integral does actually derive from the previous relation (the series expansion) of the Exponential Integral, which we write in terms of the General Exponential Integral:
$$\mathcal{E}_1(z) = \text{Ein}(z) - \ln(z) - \gamma$$
that is
$$\text{Ein}(z) = \mathcal{E}_1(z) + \ln(z) + \gamma$$
\\
where we recall that $\mathcal{E}_1(z) = -\text{Ei}(-z)$. With this relation we are immediately able to obtain an expression for $\text{Ein}(z)$ in terms of $\text{Ei}(z)$:
$$\text{Ein}(z) = - \text{Ei}(-z) + \gamma + \ln(z)$$
\\
The Modified Exponential Integral is often written as an integral, that is
$$\boxed{\text{Ein}(z) = \int_0^z \frac{1 - e^{-t}}{t}\ \de t}$$

\section{Series Expansion for $\text{Ein}(z)$}

To obtain the series expansion of the Modified Exponential Integral we start from the relation of the Generalised Exponential Integral 1: 
$$\mathcal{E}_1(z) = -\gamma - \ln(z) + \sum_{k = 1}^{+\infty}\frac{(-1)^{k+1}z^k}{k\cdot k!}$$
and making use of the above identity $\text{Ein}(z) = \mathcal{E}_1(z) + \ln(z) + \gamma$ we are able to write
$$
\begin{array}{rll}
\text{Ein}(z) & = -\gamma - \ln(z) + \displaystyle \sum_{k = 1}^{+\infty}\frac{(-1)^{k+1}z^k}{k\cdot k!} + \ln(z) + \gamma\\\\
& = \displaystyle \sum_{k = 1}^{+\infty}\frac{(-1)^{k+1}z^k}{k\cdot k!}
\end{array}
$$
The first terms of the series are:
$$\text{Ein}(z) = z - \frac{z^2}{4} + \frac{z^3}{18} - \frac{z^4}{96} + \ldots$$
With the trivial value
$$\text{Ein}(0) = 0$$
\newpage\noindent
Let's now recall the three Exponential Integral Functions we studied so far:\\\\
\textbf{\textcolor{awesome}{Exponential Integral}}
$$\text{Ei}(z) = \int_{-\infty}^z \frac{e^t}{t}\ \de t$$
\\\\
\textbf{\textcolor{awesome}{Generalised Exponential Integral 1}}
$$\mathcal{E}_1(z) = \int_z^{+\infty} \frac{e^{-t}}{t}\ \de t$$
\\
\textbf{\textcolor{awesome}{Modified Exponential Integral}}
$$\text{Ein}(z) = \int_0^z \frac{1 -e^{-t}}{t}\ \de t$$
\\
It's useful to see them plotted all together:
\\
\begin{figure}[h!]
\centering
\includegraphics[scale=1.6]{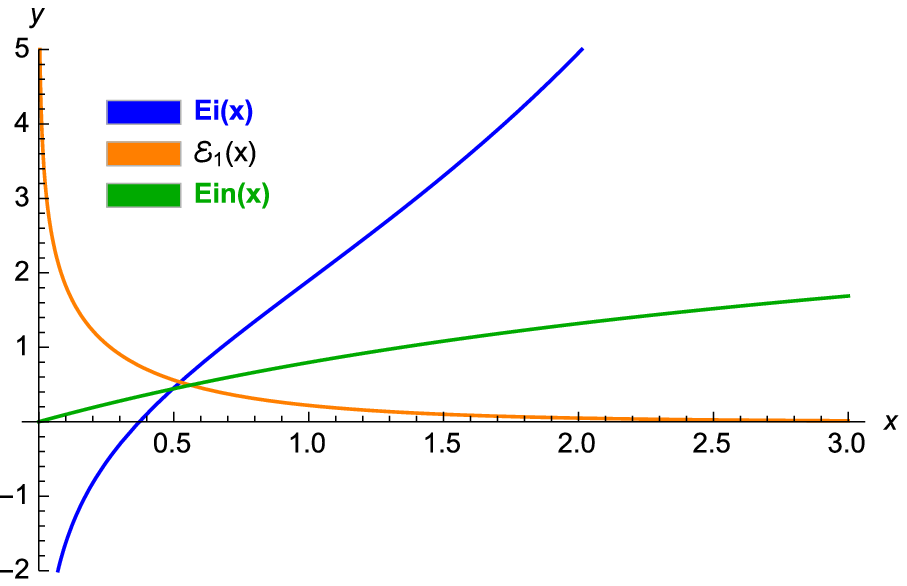}
\end{figure}
\
\\
We have to carefully notice that, contrarily to what it might appear from the plot above, \textbf{there is no common point of intersection amongst the three functions}!\\
This means \textbf{\textcolor{awesome}{there are no values for $z$ such that}}:
$$\text{Ei}(z) = \mathcal{E}_1(z) = \text{Ein}(z)$$
\newpage
\chapter{Asymptotic Series}

\section{Order Notation}

Before illustrating the asymptotic  series of the above studied function, it's quite useful to recall some definitions about the so called order notation, that is the symbols $o$, $\mathcal{O}$ and $\sim$.
\\\\
Those symbols were first used by E. Landau and P. Du Bois-Reymond and their meaning is the following: suppose $f(z)$ and $g(z)$ are functions of the continuous complex variable $z$ defined on some domain $\mathcal{D} \subset \mathbb{C}$ and possess a limit as $z\to z_0$ in $\mathcal{D}$. We hence define the following shorthand notations for the relative properties of these functions in the limit $z\to z_0$.
\\
\subsection{Asymptotically Bounded - Big O - $\mathcal{O}$}

The writing \\
$$f(z) = \mathcal{O}(g(z)) ~~~~~~~~~~~ \text{as} ~~~ z\to z_0$$
\\
means: \textit{there exist two costants $C \geq 0$ and  $\delta >0$ such that, for $0 < |z-z_0| <\delta$,}
$$|f(z)|  \leq C|g(z)|$$    
\\
We say that $f(z)$ is \textbf{asymptotically bounded} by $g(z)$ in magnitude as $z\to z_0$ or as it's usually said: "$f(z)$ is of order Big-O of $g(z)$".\\
Hence provided that $g(z)$ is not zero in a neighbourhood of $z_0$, except possibly at $z_0$, then 
\\
$$\left|\frac{f(z)}{g(z)}\right| ~~~ \text{is bounded}$$
\newpage\noindent

\subsection{Asymptotically Smaller - Little o - $o$}

The writings\\
$$f(z) = o(g(z))  ~~~~~~~~~~~ \text{as} ~~~ z\to z_0$$
\\
means: \textit{for all $\epsilon >0$ there exists $\delta > 0$ such that, for $0 < |z-z_0| < \delta$,}
$$f(z) \leq \epsilon(g(z))$$
Equivalently this means that, provided $g(z)$ is not zero in a neighbourhood of $z_0$, except possibly at $z_0$, as $z\to z_0$:
\\
$$\frac{f(z)}{g(z)}\to 0$$
We say that $f(z)$ is \textbf{asymptotically smaller} than $g(z)$ or also that $f(z)$ is of order little o of $g(z)$ as $z\to z_0$.
\\
\subsection{Asymptotically Equal - Goes Like - $\sim$}

The writings
\\
$$f(z) \sim g(z) ~~~~~~~ \text{as} ~~ z\to z_0$$
\\
means that, provided $g(z)$ is not zero in a neighbourhood of $z_0$, except possibly at $z_0$, as $z\to z_0$, 
\\
$$\frac{f(z)}{g(z)} \to 1 $$
Equivalently this means also that as $z\to z_0$
\\
$$f(z) = g(z) +  o(g(z))$$
and we say that $f(z)$ is \textbf{asymptotically equivalent} to $g(z)$ in this limit or also that $f(z)$ goes like $g(z)$ as $z\to z_0$.
\\\\
Note that the $\mathcal{O}$ order is more informative than the $o$ order about the behaviour of the function concerned as $z\to z_0$. \\
For example: 
$$\sin(z)  = z + o(z^2) ~~~~~~~ \text{as} ~~ z\to 0$$
tells us that $\sin(z)  - z$ goes to $0$ faster than $z^2$.\\
However, the writing 
$$\sin(z) = z + \mathcal{O}(z^3)$$
tells us specifically that $\sin(z) - z$ goes to zero like $z^3$. 
\newpage\noindent

\subsection{Properties of asymptotic orders}

The following properties hold, as $z\to z_0$.\\
\begin{itemize}
\item $o(f(z)) + o(f(z)) = o(f(z))$; 
\item $o(f(z)) o(g(z)) = o(f(z)g(z))$;
\item $o(o(f(z))) = o(f(z))$; \\\\
\item $\mathcal{O}(f(z)) + \mathcal{O}(f(z)) = \mathcal{O}(f(z))$; 
\item $\mathcal{O}(f(z)) \mathcal{O}(g(z)) = \mathcal{O}(f(z)g(z))$;
\item $\mathcal{O}(\mathcal{O}(f(z))) = \mathcal{O}(f(z))$; \\\\
\item $o(f(z)) + \mathcal{O}(f(z)) = \mathcal{O}(f(z))$; 
\item $o(f(z)) \mathcal{O}(g(z)) = o(f(z)g(z))$; 
\item $\mathcal{O}(o(f(z))) = o(f(z))$; 
\item $o(\mathcal{O}(f(z))) = o(f(z))$.
\end{itemize}
\
\\
\\
\subsection{Useful Examples}

\begin{itemize}
\item $f(z) = \mathcal{O}(1)$ as $z\to z_0$ means that $f(z)$ is bounded  when $z$ is close to $z_0$;
\item $f(z) = o(1)$ means that $f(z) \to 0$ as $z\to z_0$, namely that $f(z)$ is infinitesimal; 
\item $\ln(z) = \mathcal{O}(z^{-a})$ as $z\to 0^+$ for $a >0$;
\item $\ln(x) = o(z^a)$ as $z\to +\infty$; 
\item $\sin(z) \sim z$ as $z\to 0$;
\item $\sin(z) = \mathcal{O}(1)$ for every $z\in\mathbb{R}$;
\item as $z\to 0 ~~~~~ \begin{cases} \sin\left(\frac{1}{z}\right)  = \mathcal{O}(1) \\ \cos(z)  \sim 1 - \frac{1}{2}z^2 \end{cases}$ \\
\item as $z\to 0^+ ~~~~~ \begin{cases} z^2 = o(t) \\ e^{-1/z} = o(1) \\ \tan(z) =  \mathcal{O}(z) \\ \sin(z) \sim z \end{cases}$\\
\item as $z\to +\infty ~~~~~ \begin{cases} z^{1000} = o(e^t) \\ \cos(z) = \mathcal{O}(1) \end{cases}$\\
\item If we have a polynomial function like: $f(z) = 5z^2  + z + 3$ then: \\\\
 as $z\to +\infty  ~~~~~ \begin{cases} 
 f(z) = o(z^3) \\
 f(z) = \mathcal{O}(z^2) \\
 f(z) \sim 5z^2
 \end{cases}$
\end{itemize}
\
\\
\subsection{Asymptotic Sequences}

Let $M$ be a set of real or complex numbers, with a limit point $z_0$. \\
Let $\phi_n : M \to \mathbb{R}$ or also $\mathbb{C}$, and $n\in\mathbb{N}$.  Let it also be that $\phi_n(z) \neq 0$ in a neighbourhood $I_n$ of $z_0$. \\\\
The sequence $\{\phi_n\}$ is called \textbf{asymptotic sequence} at $z\to z_0$ if $\forall n\in\mathbb{N}$ \\
$$\phi_{n+1}(z) = o(\phi_n(z))$$
\\
\subsection{Asymptotic Series}

Let $f : M\to \mathbb{R}$ and $z_0$ be a limit point in $M$.\\
Let $\{\phi_n\}$ be an asymptotic sequence as $z\to z_0$, $z\in M$.\\
We say that the function $f$ is \textbf{expanded in an asymptotic series}
$$f(z) \sim \sum_{n = 0}^{+\infty} \alpha_n \phi_n(z) ~~~~~~~~~~~~~ \alpha_n = \text{constants}$$
if $\forall N\geq 0$
$$R_N(z) \equiv f(z)  -   \sum_{n = 0}^N \alpha_n \phi_n(z) = o(\phi_N(z))$$
\\
The series is called the \textbf{asymptotic expansion} of the function $f$ with respect to the asymptotic sequence $\{\phi_n\}$.\\
$R_N(z)$ is called the remainder or the rest term of the asymptotic series.
\\\\\\
\textbf{Remarks.}
\\
\begin{enumerate}
\item The condition $R_N(z) = o(\phi_N(z))$ means, in particular, that
$$\lim_{z\to z_0} R_N(z) = 0 ~~~~~~~ \text{for any fixed} ~~  N$$
\item \textcolor{awesome}{The asymptotic series could diverge}. This happens if
$$\lim_{N\to +\infty} R_N(z) \neq 0 ~~~~~~~ \text{for some fixed} ~~  z$$
\item There are, in general, three possibilities:
\begin{itemize}
\item The asymptotic series converges to $f(z)$; 
\item The asymptotic series converges to $g(z) \neq f(z)$; 
\item The asymptotic series diverges.
\end{itemize}
\end{enumerate}
\
\\
It's worthwhile to recall a very important \textcolor{azure}{theorem} and a remark:\\\\
\textit{The asymptotic expansion of a function with respect to an asymptotic sequence is unique.}
\\\\
\textit{Two different functions can have the same asymptotic expansion.}
\\\\
To see an example of the last remark, let's consider the functions
$$f(z) = e^z ~~~~~~~~~~~~~ g(z) = e^{z} + e^{-1/z}$$
With respect to the asymptotic sequence $\{z^n\}$ they have the same asymptotic expansion:
$$e^z \sim e^z + e^{-1/z} \sim \sum_{n = 0}^{+\infty} \frac{z^n}{n!} ~~~~~~~ z\to 0^+$$
Although this was just an introduction to asymptotic series, we are ready to face the asymptotic series of the Exponential Integral.

\section{Asymptotic Series of $\text{Ei}(z)$ and $\mathcal{E}_1(z)$}

To obtain the asymptotic series, we just integrate $\text{Ei}(z)$ by parts $n$ times:
\begin{equation*}
\begin{split}\int_{-\infty}^z \frac{e^t}{t}\ \de t & = \int_{-\infty}^z \frac{1}{t}\ \de\left(e^t\right) \\\\
& = \frac{e^z}{z} + \int_{-\infty}^z \frac{e^t}{t^2}\ \de t \\\\
& = \frac{e^z}{z} + \frac{e^z}{z^2} + 2\int_{-\infty}^z \frac{e^t}{t^3}\ \de t \\\\
& = \ldots \\\\
& = e^z\left\{\frac{1}{z} + \frac{1}{z^2} + \frac{1\cdot 2}{z^3} + \ldots + \frac{n!}{z^{n+1}}\right\} + (n+1)!\int_{-\infty}^z \frac{e^t}{t^{n+2}}\ \de t
\end{split}
\end{equation*}
\
\\
From which it follows the \textbf{\textcolor{awesome}{Asymptotic Series for the Exponential Integral}}
$$\boxed{\text{Ei}(z) = \frac{e^z}{z}\left[\sum_{k = 0}^{n}\frac{k!}{z^k} + \Upsilon_n(z)\right]}$$
where
$$\Upsilon_n(z) = (n+1)!z e^{-z}\int_{-\infty}^z \frac{e^t}{t^{n+2}}\ \de t$$
is the remainder.
\\
\\
We can appraise the remainder as an order (meant as the Big-$\mathcal{O}$ notation); in particular we find
$$|\Upsilon_n(z)| = \mathcal{O}\left(|z|^{-n-1}\right)$$
And we can then write
$$\text{Ei}(z) = \frac{e^z}{z}\left(\sum_{k = 0}^n \frac{k!}{z^k} + \mathcal{O}(|z|^{-n-1})\right)$$
\\\\
By the relationship we wrote many times between $\text{Ei}(z)$ and $\mathcal{E}_1(z)$, namely $\mathcal{E}_1(z) = - \text{Ei}(-z)$, we are able to obtain also the \textbf{\textcolor{awesome}{asymptotic series of $\mathcal{E}_1(z)$}}.\\
$$\mathcal{E}_1(z) = -\frac{e^{-z}}{-z}\left(\sum_{k = 0}^n \frac{k!}{(-z)^k} + \mathcal{O}(|z|^{-n-1})\right)$$
That is
$$\boxed{\mathcal{E}_1(z) = \frac{e^{-z}}{z}\left(\sum_{k = 0}^n \frac{(-1)^k k!}{z^k} + \mathcal{O}(|z|^{-n-1})\right)}$$
\\
We can give a more rigorous proof by starting from its own definition as an integral, integrating by parts $n$ times and obtaining the same result. We will do this also to show some particular noteworthy questions about the partial sum of the Exponential Integral and the rest term.
\\\\
Starting from $$\mathcal{E}_1(z) = \int_z^{+\infty} \frac{e^{-t}}{t}\ \de t$$
we integrate by parts $n$ times, to get a series expansion (for $z >> 1$):
\begin{equation*}
\begin{split}
\mathcal{E}_1(z) & = \frac{e^{-z}}{z} - \frac{e^{-z}}{z^2} + 2\int_z^{+\infty} \frac{e^{-t}}{t^3}\ \de t \\\\
& = e^{-z}\underbrace{\left[\frac{1}{z} - \frac{1}{z^2} + \ldots  + (-1)^{n-1} \frac{(n-1)!}{z^n}\right]}_{S_n(z)} + \underbrace{(-1)^n n! \int_z^{+\infty} \frac{e^{-t}}{t^{n+1}}\ \de t}_{\Upsilon_n(z)}
\end{split}
\end{equation*}
\
\\
Where we denoted as $S_n(z)$ the partial sum and as $\Upsilon_n(z)$ the rest term.\\\\
It's important to note that $S_n(z)$ is divergent for every $z$.\\
$\Upsilon_n(z)$ is also unbounded as $n\to +\infty$, but the sum $S_n(z) + \Upsilon_n(z)$ must be bounded since $\mathcal{E}_1(z)$ is bounded and definite for every $z>0$.
\\\\
Let's take a look at the rest term.\\
Be $n$ fixed and let $z$ become larger:
\begin{equation*}
\begin{split}
|\Upsilon_n(z)| & = \left|(-1)^n n! \int_z^{+\infty} \frac{e^{-t}}{t^{n+1}}\ \de t\right| \\\\
& = |(-1)^n n| \int_z^{+\infty} \frac{e^{-t}}{t^{n+1}}\ \de t \\\\
& = n! \int_z^{+\infty} \frac{e^{-t}}{t^{n+1}}\ \de t  \\\\
& < \frac{n!}{z^{n+1}} \int_z^{+\infty} e^{-t}\ \de t\\\\
& = \frac{n!}{z^{n+1}} e^{-z}
\end{split}
\end{equation*}
which tends to $0$ rapidly as $z\to +\infty$. \newpage \noindent 
Also note that
\begin{equation*}
\begin{split}
\left|\frac{\Upsilon_n(z)}{(n-1)! e^{-z} z^n}\right| & = \frac{|\Upsilon_n(z)|}{(n-1)! e^{-z} z^n} \\\\
& < \frac{n! e^{-z}z^{-(n+1)}}{(n-1)! e^{-z}z^{-n}} \\\\
& = \frac{n}{z}
\end{split}
\end{equation*}
which also goes to zero as $z\to +\infty$.\\
Thus 
$$\mathcal{E}_1(z) = S_n(z) + o(\text{last term in}~~S_n(z)) ~~~~~~~ \text{as} ~~ z\to +\infty$$
\\
For large $z$ and $n$ fixed, $S_n(z)$ is a very good approximation to $\mathcal{E}_1(z)$; the accuracy of the approximation increases as $z$ increases.
$$\mathcal{E}_1(z) \sim e^{-z}\left(\frac{1}{z} - \frac{1}{z^2} + \frac{2!}{z^3} + \ldots \right)$$
That is
$$\mathcal{E}_1(z) = \frac{e^{-z}}{z}\sum_{k = 0}^n \frac{(-1)^k k!}{z^k}$$
as obtained before.\\
To end this section, we show the plot of the asymptotic series compared to the integral definition:
\\
\begin{figure}[h!]
\centering
\includegraphics[scale=1.5]{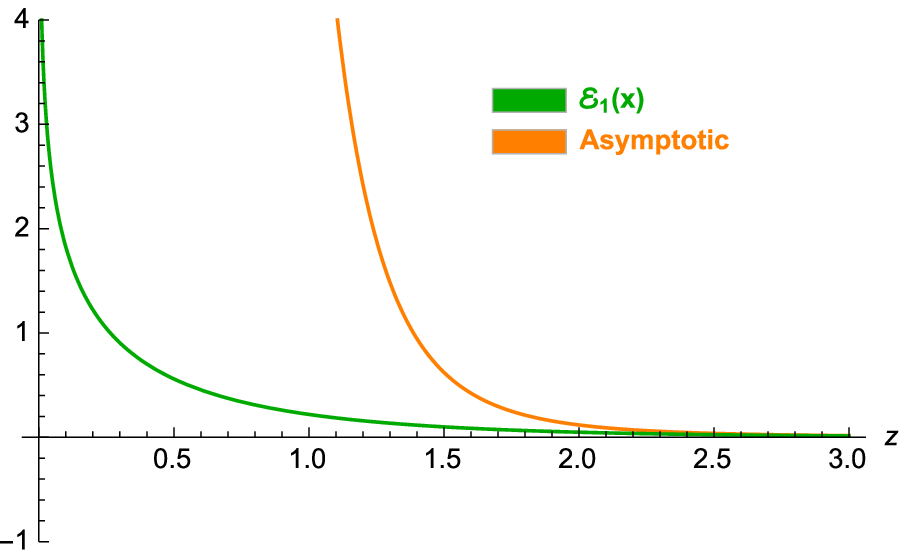}
\end{figure}
\
\\

\section{Imaginary Argument}

Let's now face a very interesting particular case, that is when the argument of the Exponential Integral (we will refer hereafter to $\mathcal{E}_1(z)$, unless spoken differently) is purely imaginary, \textit{i.e.} $z = ix$. Starting from the definition
$$\mathcal{E}_1(z) = \int_z^{+\infty} \frac{e^{-t}}{t}\ \de t$$
We firstly redefine the variable of integration as
$$t\to tz ~~~~~~~ \de t = z\ \de t$$
Obtaining
$$\mathcal{E}_1(z) = \int_1^{+\infty} \frac{e^{-tz}}{t}\ \de t$$
At this point, setting $z = ix$ we easily find:
$$\mathcal{E}_1(ix) = \int_1^{+\infty}\frac{-i\sin(tx) + \cos(tx)}{t}\ \de t = \int_1^{+\infty}-i \frac{\sin(tx)}{t}\ \de t + \int_1^{+\infty} \frac{\cos(tx)}{t}\ \de t$$
\\
The first term can be manipulated a bit to get a well known function which we foretell to be the so called \textbf{Integral Sine} function.
$$
-i\int_1^{+\infty} \frac{\sin(tx)}{t}\ \de t ~~~~~~~~~ 
\begin{cases}
t \to \frac{t}{x} \\\\
\de t = \frac{1}{x}\ \de t
\end{cases}
~~~~~ \longrightarrow ~~~~~ -i\int_x^{+\infty} \frac{\sin t}{t}\ \de t
$$
The integral from $x$ to $+\infty$ can be rewritten as the integral over the whole $\mathbb{R}^+$ range minus the missing range $[0, x]$ as follows:
$$\int_x^{+\infty} = \int_0^{+\infty} - \int_0^x$$
Obtaining
$$-i\left[\int_0^{+\infty}\frac{\sin t}{t}\ \de t - \int_0^x \frac{\sin t}{t} \right]$$
The very first term is a well known integral, sometimes called \textit{the Dirichlet Integral}, which can be easily computed with the help of Residues Calculus:
$$\int_0^{+\infty}\frac{\sin t}{t}\ \de t = \frac{\pi}{2}$$
whereas the second term is the function we foretold before to occur, that is the \textbf{\textcolor{awesome}{Sine Integral Function}}
$$\boxed{\text{Si}(x) = \int_0^x\frac{\sin t}{t}\ \de t}$$
So we ended up with
$$-i \int_1^{+\infty}\frac{\sin t}{t}\ \de t = -i\frac{\pi}{2} + i\text{Si}(x)$$
\\\\
For what concerns the second term in the $\mathcal{E}_1(ix)$ we can reason in a quite similar way:
$$\int_1^{+\infty} \frac{\cos t}{t}\ \de t ~~~
\begin{cases}
t \to \frac{t}{x} \\\\
\de t = \frac{1}{x}\ \de t
\end{cases}
~~~~~ \longrightarrow ~~~~~ 
\int_x^{+\infty} \frac{\cos t}{t}\ \de t
$$
We cannot use the same trick we used before, because this last integral does not converge in the range $[0, +\infty]$. The form we obtained is what it's called the \textbf{\textcolor{awesome}{Cosine Integral Function}}
$$\boxed{\text{Ci}(x) = -\int_x^{+\infty} \frac{\cos t}{t}\ \de t}$$
Usually written as
$$\boxed{\text{Ci}(x) = \int_{+\infty}^x \frac{\cos t}{t}\ \de t}$$
\\ 
\textcolor{azure}{\hrule}
\noindent \ \\
\textbf{Remark}: For the sake of completeness, it's worthwhile to mention that there are two different Sine Integral Functions, that is:
$$\text{Si}(x) = \int_0^x\frac{\sin t}{t}\ \de t$$
$$\text{si}(x) = - \int_x^{+\infty}\frac{\sin t}{t}\ \de t$$
The second one of which will not be treated in this paper.
\\
\textcolor{azure}{\hrule}
\noindent 
\
\\
Due to what we have done so far, we are able to write the Exponential Integral with an imaginary argument in terms of the Sine Integral and the Cosine Integral functions:
$$\mathcal{E}_1(ix) = -i\frac{\pi}{2} + i\text{Si}(x) - \text{Ci}(x)$$
It's very important and cool to notice that despite the imaginary argument, the Sine and Cosine Integral functions are purely real!
\\\\
\textcolor{azure}{\hrule}
\noindent 
\
\\
\textbf{Remark}: It may be useful to know that the Cosine Integral Function is often expressed as:
$$\text{Ci}(x) = \gamma + \ln(x) + \int_0^x \frac{\cos t - 1}{t}\ \de t$$
where $\gamma$ is the Euler-Mascheroni constant.
\\
\textcolor{azure}{\hrule}
\noindent 
\newpage 
\noindent 
As usual, it's important to visualize the plots of those two special functions, the first one of which we already met during the treatise of the Gibbs Phenomenon (\textit{cfr:} Enrico Masina - Note sulle Serie di Fourier).
\\
\begin{figure}[h!]
\centering
\includegraphics[scale=1.7]{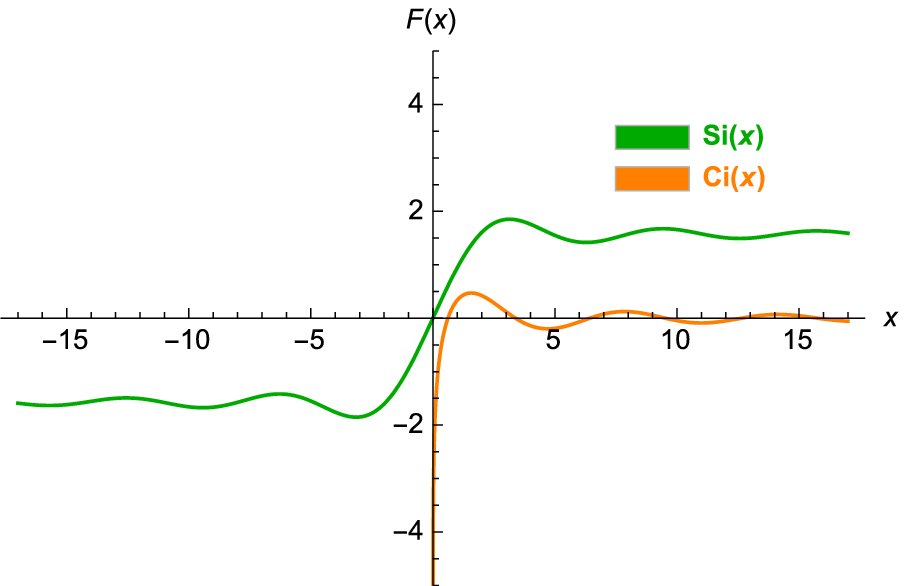}
\end{figure}
\
\\
Here are some notable properties and values:
$$\text{Si}(-z) = -\text{Si}(|z|) ~~~~~~~ \text{Ci}(-z) = \text{Ci}(z) - i\pi$$
$$\text{Si}(\infty) = \frac{\pi}{2} ~~~~~~~ \text{Si}(0) = 0$$
$$\text{Ci}(\infty) = 0 ~~~~~~~~~~~ \text{Ci}(0^+) = -\infty$$
\\\\
What we have done so far for $\mathcal{E}_1(z)$ can obviously be done in the same way for the classic Exponential Integral $\text{Ei}(z)$. We give a quick input on how to get the result.
\\\\
Directly substituting $z = ix$ into the definition, one gets:
$$\text{Ei}(ix) = \int_{-\infty}^{z} \frac{e^t}{t}\ \de t =  -\int_{ix}^{+i\infty} \frac{e^t}{t}\ \de t$$
The trick now is to perform an imaginary change of variable, which is very similar to what is called a "Wick Rotation" in theoretical physics:
$$t\to i u ~~~~~~~ \de t = i\ \de u$$
Which leads us to
\begin{equation*}
\begin{split}
\text{Ei}(ix) & = \int_{+\infty}^x \frac{e^{iu}}{u}\ \de u \\\\
& = \int_{+\infty}^x \frac{\cos u}{u}\ \de u + i\int_{+\infty}^x \frac{\sin u}{u}\ \de u
\end{split}
\end{equation*}
\ \\
as wanted.
\\\\
\subsubsection{\textcolor{awesome}{Some important relations.}}

We present a short list of the most important and useful relations between the Sine Integral, the Cosine Integral, $\text{Ei}$ and $\mathcal{E}_1(x)$.\\
$$\text{Ei}(ix) = \text{Ci}(x) - i\left(\frac{\pi}{2} - \text{Si}(x)\right) ~~~~~~~ x > 0$$
$$\text{Ei}(-ix) = \text{Ci}(x) + i\left(\frac{\pi}{2} - \text{Si}(x)\right) ~~~~~~~ x > 0$$
\\
$$\text{Si}(ix) = \frac{i}{2}\left[\text{Ei}(x) + \mathcal{E}_1(x)\right] ~~~~~~~ x>0$$
$$\text{Ci}(ix) = \frac{1}{2}\left[\text{Ei}(x) - \mathcal{E}_1(x)\right] +i\frac{\pi}{2} ~~~~ x>0$$

\subsection{Series Expansion of $\text{Si}(x)$ and $\text{Ci}(x)$}

It might be interesting to show the series expansions for these two special functions, without giving the direct proof (which is however quite straightforward to obtain).
\\ \\
$$\text{Si}(z) = \sum_{k = 0}^{+\infty} \frac{(-1)^k z^{2k+1}}{(2k+1)(2k+1)!}$$
\\
$$\text{Ci}(z) = \gamma + \ln(z) + \sum_{k = 0}^{+\infty} \frac{(-1)^k z^{2k}}{(2k)\ (2k)!}$$
\\\\
These series are convergent at any complex $z$, although for $|z| >> 1$ the series will initially converge slowly, which means more terms are required for a high precision.
\newpage\noindent 
\subsection{Asymptotic Series for $\text{Si}(x)$ and $\text{Ci}(x)$}

The following series are asymptotic and divergent, although they can be used for an estimate or even a precise evaluation at $\Re (x) >> 1$:
\\
$$\text{Si}(z) = \frac{\pi}{2} - \frac{\cos z}{z}\left(\sum_{k = 0}^{+\infty} (-1)^k\frac{(2k)!}{z^{2k}}\right) - \frac{\sin z}{z} \left(\sum_{k = 0}^{+\infty}(-1)^k \frac{(2k+1)!}{z^{2k+1}}\right)$$
\\
$$\text{Ci}(z) = \frac{\sin z}{z} \left(\sum_{k = 0}^{+\infty}(-1)^k \frac{(2k+1)!}{z^{2k+1}}\right) - \frac{\cos z}{z}\left(\sum_{k = 0}^{+\infty} (-1)^k\frac{(2k)!}{z^{2k}}\right)$$

\chapter{The Logarithmic Integral}

A very important function connected to the Exponential Integral function is the Logarithmic Integral function, also known as Integral Logarithm, denoted by $\text{li}(x)$.\\
This function is extremely important and relevant in a lot of problems in physics and it also has a large importance in the mathematical branch of Number Theory, occurring for example in the prime number theorem: it's indeed used to estimate  the number of prime numbers less than a given value.
\\\\
The Logarithmic Integral has an integral representation defined for all positive real numbers $x \neq 1$ by the definite integral: 
$$\boxed{\text{li}(x) = \int_0^x \frac{\de t}{\ln(t)}}$$
\\
Because of the singularity at $t = 1$, we must interpret this integral as a  Cauchy Principal Value:
$$\text{li}(x) = \text{P.V.}\int_0^x\frac{\de t}{\ln(t)} = \lim_{\alpha\to 0^+}\left(\int_0^{1-\alpha} \frac{\de t}{\ln(t)} + \int_{1+\alpha}^x \frac{\de t}{\ln(t)}\right)$$
\\
A definition for the Logarithmic Integral with complex argument can however be given as
$$\text{li}(z) = \int_0^z \frac{\de t}{\ln(t)} 
~~~~~~~~
\begin{cases}
|\arg(z)| < \pi \\\\
|\arg(1-z)| < \pi 
\end{cases}
$$
\\
The connection with the Exponential Integral can be seen by performing this change of variables:
$$u = \ln(t) ~~~~~~~ \de u = \frac{\de t}{t}$$
leading us to
$$\text{li}(z) = \int_{-\infty}^{\ln(z)} \frac{e^u}{u}\ \de u$$
Where we recognize the integral to be exactly the Exponential Integral $\text{Ei}(z)$ function of argument $\ln(z)$, \textit{viz.}
$$\boxed{\text{li}(z) = \text{Ei}(\ln(z))}$$
\\
As always, we plot the function to have a clearer idea of its behaviour.
\\
\begin{figure}[h!]
\centering
\includegraphics[scale=1.7]{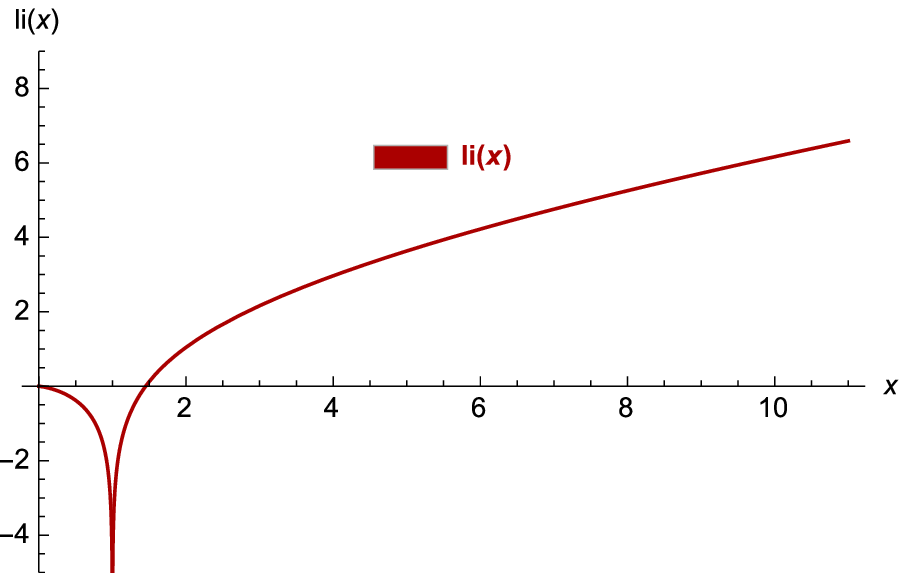}
\end{figure}
\
\\
We can easily notice the singularity at $x = 1$.\\\\
The function $\text{li}(z)$ has a single positive real zero which occurs at 
$$x = 1.4513692348(\ldots)$$
which is known to be the \textit{Ramanujan-Soldner constant}.
\\
\section{Series of $\text{li}(z)$}

The Series expansion of $\text{li}(z)$ is easily deducible by putting $\ln(z)$ 
as the argument in the series of $\text{Ei}(z)$, \textit{viz.}
\\
$$\text{li}(z) = \gamma + \ln(-\ln(z)) + \sum_{k = 1}^{+\infty} \frac{(\ln(z))^k}{k\cdot k!}$$
There is a more rapidly convergent series, due to S. Ramanujan which reads:
\\
$$\text{li}(z) = \gamma +\ln(\ln(z))  + \sqrt{z}\sum_{k = 1}^{+\infty} \frac{(-1)^{k-1} (\ln(z))^k}{k! 2^{k-1}} \sum_{j = 0}^{\floor*{\frac{k-1}{2}}} \frac{1}{2j+1}$$
Where the notation $\floor*{\frac{k-1}{2}}$ represents "the Floor function".
\newpage\noindent 
\section{Modified Logarithmic Integral}

If we take the function $\mathcal{E}_1(z)$ with argument $z\to \ln(z)$ we obtain the so called \textbf{\textcolor{awesome}{Modified Logarithmic Integral}}: 
\\
$$\text{li}_1(z) = \mathcal{E}_1(\ln(z)) = -\text{Ei}(-\ln(z))$$
\\
So it's actually not a big deal, but it was worth mentioning as it shows how all the Exponential Integral functions are related.

\section{Asymptotic expansion of $\text{li}(z)$ and $\text{li}_1(z)$}

Using the previous relation between the Logarithmic Integral and the Exponential Integral we are easily able to show its asymptotic series.
\\
$$\text{li}(z) \sim \frac{z}{\ln(z)}\left[\sum_{k = 0}^n \frac{k!}{(\ln(z))^k} + \Upsilon_n(z)\right]$$
\\
Where the estimated remainder is
$$\Upsilon_n(z) = \mathcal{O}\left(|\ln(z)|^{-n-1}\right)$$
\\\\
In the same way, using $\mathcal{E}_1(\ln(z))$ we immediately find the asymptotic series for $\text{li}_1(z)$: 
\\
$$\text{li}_1(z) \sim \frac{1}{z\ln(z)}\left[\sum_{k = 0}^{+\infty} \frac{(-1)^k k!}{(\ln(z))^k} + \Upsilon_n(z)\right]$$

\section{Offset Logarithmic Integral}

We are finally arrived to one of the most important special functions of this paper, which is the so called \textbf{Offset Logarithmic Integral} or sometimes called \textbf{Eulerian Logarithmic Integral}, which is defined (we now mind only about real argument) as\\
$$\text{Li}(x) = \text{li}(x) - \text{li}(2)$$\\
or via integral representation as\\
$$\text{Li}(x) = \int_2^{x} \frac{\de t}{\ln(t)}$$\\
The advantage of such a representation is that the Offset L.I. avoids the singularity in the domain of integration. \\\\
For the records
$$\text{li}(2) \approx 1.045163780117492784 844588 889194 613136 522615 578151$$
\\
This function is extremely important in the field of Number Theory, and in particular way in the theory of distribution of prime numbers; this function represents indeed a very good estimation for the numbers of prime number before a given value, and it is also one of the best approximations for the (still unknown) \textbf{prime numbers counter function} $\pi(x)$.
\\
Indeed since the exact form of the function $\pi(x)$ is unknown (provided it have an exact form), this function is a great approximation also because as $x\to +\infty$ its asymptotic behaviour is defined as
$$\text{Li}(x) = \mathcal{O}\left(\frac{x}{\ln(x)}\right)$$
\\
It's quite obvious to imagine what the plot looks like, since it's basically the Logarithmic integral with a translation of a quantity equal to $\text{li}(2)$ with respect to the \textsf{Y} axis, in order to set
$$\text{Li}(2) = 0$$
\\
\begin{figure}[h!]
\centering
\includegraphics[scale=1.5]{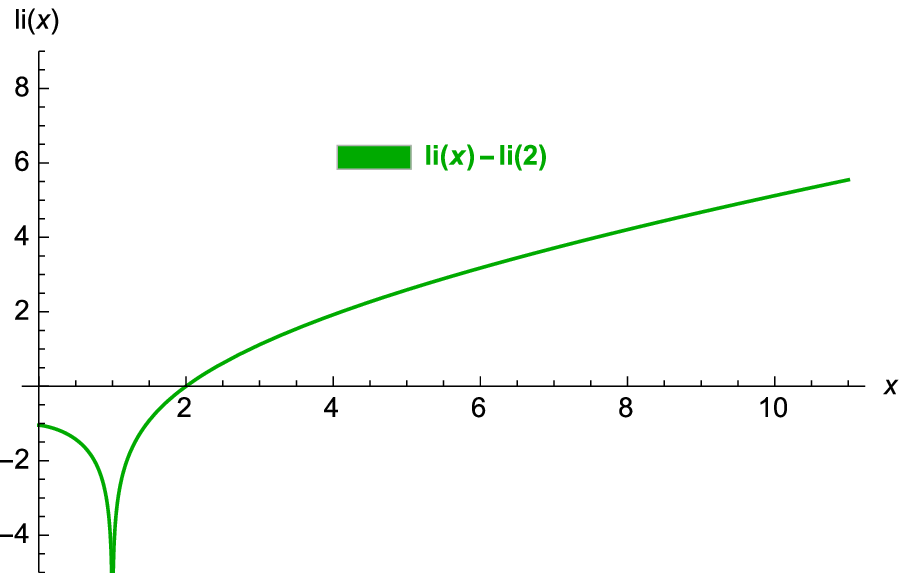}
\end{figure}
\
\\
The comparison with the plot at page $27$ shows the slight but important difference between the two functions. 

\section{Prime Numbers counter $\pi(x)$ and $\text{Li}(x)$}

We spoke about its role in the distribution of prime numbers, and in particular about its being a very good estimation for the $\pi(x)$ function; in $1976$ Lowell Schonfeld showed that the error in choosing $\text{Li}(x)$ as a good candidate for $\pi(x)$ is:
\\
$$|\pi(x) - \text{Li}(x)| < \frac{\sqrt{x}\ln(x)}{8\pi} ~~~~~~~ \text{for} ~~ x > 2657$$
It may appear a rough estimation, but let's make a numerical example: let's set $x = 10000$.\\
Thanks to numerical tables, we know that the number of prime numbers less than $10000$ is
$$\pi(10000) = 1229$$
hence we have the exact value. But we don't know the form of the function $\pi(x)$ and if we use the Offset L.I. we get:
$$\text{Li}(10000) = \int_2^{10000}\frac{\de t}{\ln(t)} = 1245$$
(the correct result would be $1254.09$).\\\\
Hence the error is $16$.
\\
Comparing with the Schonfeld estimation we find
\begin{equation*} 
\begin{split}
|\pi(x) - \text{Li}(x)| & < \frac{\sqrt{10000}\ln(10000)}{8\pi} \\\\
& < 36.64677
\end{split}
\end{equation*}
which is true.
\\\\
We present now the plot of the two functions, $\pi(x)$ and $\text{Li}(x)$, to show how good is the approximation (although it is not the best one!).
\\\\
\textbf{Plot $\pi(x)$, $\text{Li}(x)$ for $x = 100$}\\\\
\begin{figure}[h!]
\centering
\includegraphics[scale=1.4]{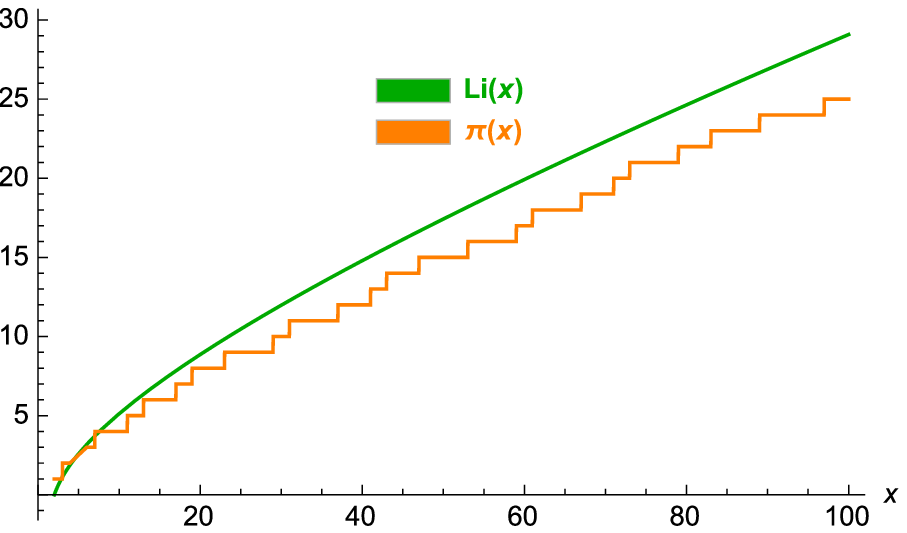}
\end{figure}
\ \\
\newpage \noindent

\textbf{Plot $\pi(x)$, $\text{Li}(x)$ for $x = 200$}
\\\\
\begin{figure}[h!]
\centering
\includegraphics[scale=1.4]{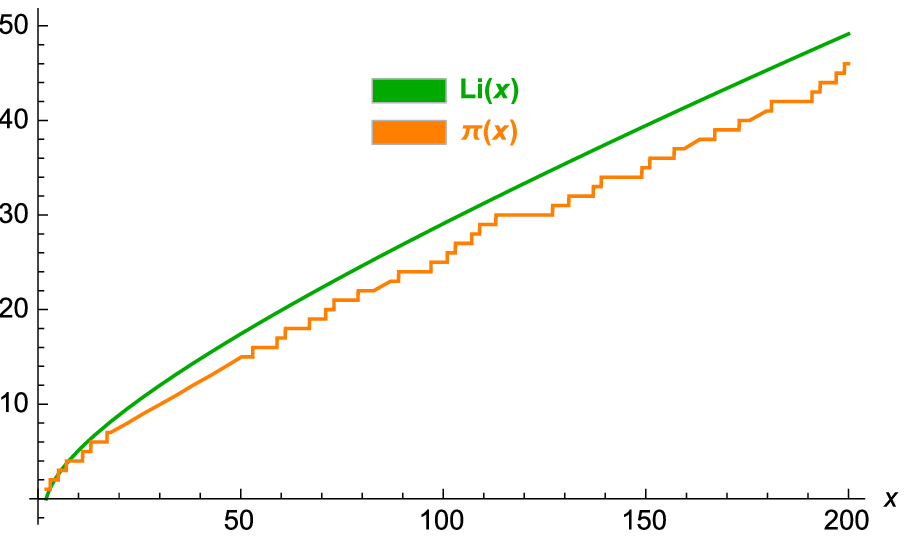}
\end{figure}
\
\\ \\
\textbf{Plot $\pi(x)$, $\text{Li}(x)$ for $x = 400$}\\\\
\begin{figure}[h!]
\centering
\includegraphics[scale=1.4]{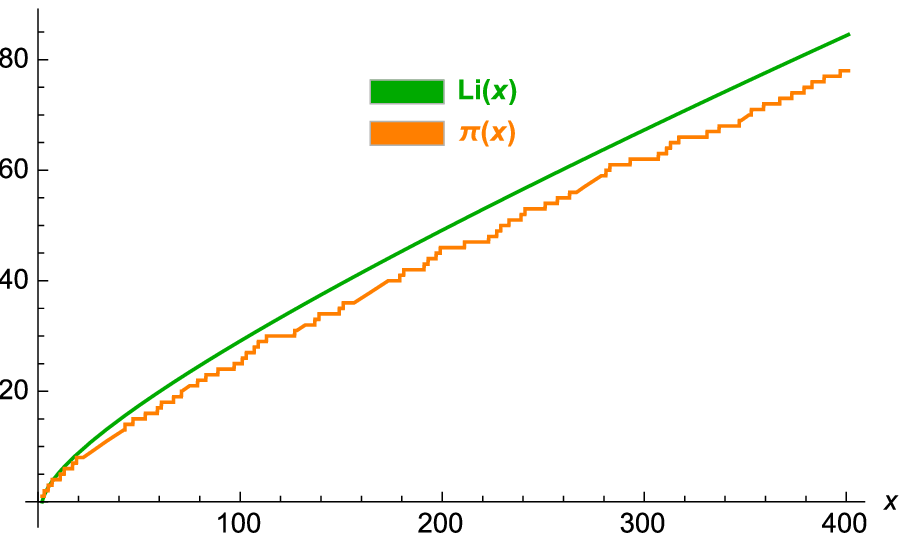}
\end{figure}
\
\\
\newpage \noindent
\textbf{Plot $\pi(x)$, $\text{li}(x)$ for $x = 700$}\\\\
\begin{figure}[h!]
\centering
\includegraphics[scale=1.4]{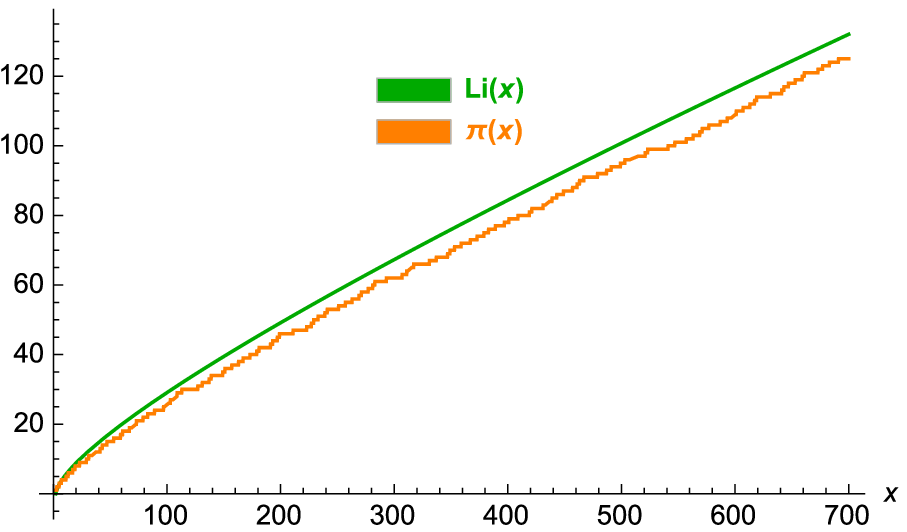}
\end{figure}
\
\\ \\
Another estimation for the Offset Logarithmic Function is given by
\\
$$\frac{\text{Li}(x)}{x/\ln(x)} \sim 1 + \frac{1}{\ln(x)} + \frac{2}{(\ln(x))^2} + \frac{6}{(\ln(x))^3} + \ldots \sim \sum_{k = 0}^{+\infty} \frac{k!}{(\ln(x))^k} $$
which gives a more accurate asymptotic behaviour:
$$\text{Li}(x) - \frac{x}{\ln(x)} = \mathcal{O}\left(\frac{x}{(\ln(x))^2}\right)$$
Other estimations for $\pi(x)$ have been computed in the past, years before Schonfeld's one, and it's a good thing to mention them:
$$\pi(x) \sim \text{Li}(x) + \mathcal{O}\left(\sqrt{x}\ln(x)\right) ~~~~~~~~~~~ \text{1901 - Von Koch}$$
$$\pi(x) \sim \text{Li}(x) + \mathcal{O}\left(xe^{-a\sqrt{\ln(x)}}\right) ~~~~~~~~~~~ \text{1899 - La Vallé Poussin}$$
\\\\
Besides the Prime Number topic, which is a huge and still unexplored field alone, there would be many more topics to discuss about in which the Exponential Integral and the Logarithmic Integral functions (and related special functions) take part in a very active and important way,  but it will be for another time.

\chapter{Physical Application and Exercises}

\section{Application}

As a simple example of an application of the special functions studied in this paper, we consider the electromagnetic energy radiated by a linear oscillator of length $2\ell = \frac{\lambda}{2}$, driven by an alternating current $I$ of frequency $\omega = \frac{2\pi c}{\lambda}$ (where $c$ is the velocity of light and $\lambda$ is the wavelenght), whose distribution along the conductor is given by
$$I = I_0\cos \left(\frac{\pi z}{2\ell}\right)\cos(\omega t) ~~~~~~~ -\ell \leq z \leq \ell$$
\\
Let $\mathbf{E}(t)$ and $\mathbf{H}(t)$ denote the time-dependent electric and magnetic field vectors, with complex amplitudes $\mathbf{E}$ and $\mathbf{H}$, so that\\
$$\mathbf{E}(t) = \Re\{\mathbf{E}e^{i\omega t}\}, ~~~~~~~ \mathbf{H}(t) = \Re\{\mathbf{H}e^{i\omega t}\}$$
\\
Then the power radiated by the oscillator, averaged over a period $T = \lambda/c$, is given by the well known formula
$$\mathcal{P} = \Re\left\{\frac{c}{8\pi} \int_S (\mathbf{E}\times\mathbf{H}^*)\cdot \mathbf{n}\ \de S\right\}$$
\\
where $S$ is an arbitrary surface surrounding the oscillator, $\mathbf{n}$ is the exterior normal to $S$, and $\mathbf{H}^*$ is the vector whose components are the complex conjugates of those of $\mathbf{H}$.
\\\\
In the present case, the vectors $\mathbf{E}$ and $\mathbf{H}$ have components $(E_r, E_{\theta}, 0)$, and $(0, 0, H)$, in a spherical coordinates system $(r, \theta, \phi)$.\\
For $S$ it's convenient to choose a sphere or arbitrarily large radius $r = \rho$. Then we get
$$\mathcal{P} = \Re\left\{\frac{c\rho^2}{4}\int_0^{\pi} E_{\theta} H^*\sin\theta\ \de \theta\right\}$$
where $H^*$ is the complex conjugate of $H$.\\
We can replace the exact values of $E_{\theta}$ and $H$ by their asymptotic expressions for large $r$. Using the well known formulas for the components of the electromagnetic field of an elementary dipole (\textit{see as a ref. :} G. Joos - Theoretical Physics, pp. 338-340), and integrating with respect to $z$ we easily find
$$H \approx E_{\theta} \approx \frac{I_0 i k}{c\rho}e^{-ik\rho}\sin\theta \int_{-\ell}^{+\ell}\cos \frac{\pi z}{2\ell} e^{ikz\cos\theta}\ \de z$$
which is equal to
$$H \approx E_{\theta} \approx \frac{2I_0 i}{c\rho} e^{-ik\rho}\frac{\cos\left(\frac{1}{2}\pi\cos\theta\right)}{\sin\theta}$$
\\
For sufficiently large $\rho$.\\
$k = \omega/c$.\\\\
It follows that
$$\mathcal{P} = \frac{I_0^2}{c}\int_0^{\pi} \frac{\cos^2\left(\frac{1}{2}\pi\cos\theta\right)}{\sin\theta}\ \de\theta$$
in which we made use of $\ell = \frac{\lambda}{4} = \frac{\pi c}{2\omega} = \frac{\pi}{2k}$.
\\\\
The integral can be expressed in terms of the Cosine Integral $\text{Ci}(x)$. In fact, introducing the new variable of integration $x = \cos\theta$ we have
\\
$$\mathcal{P} = \frac{I_0^2}{c}\int_0^{1} \frac{1 + \cos \pi x}{1 - x^2}\ \de x = \frac{I_0^2}{c}\left(\int_0^{1} \frac{1 + \cos \pi x}{1 - x}\ \de x + \int_0^{1} \frac{1 + \cos \pi x}{1 + x}\ \de x \right)$$
Introducing $y = 1-x$ we can rewrite it all as
$$\mathcal{P} = \frac{I_0^2}{c}\left(\int_0^{1} \frac{1 + \cos \pi y}{y}\ \de y + \int_1^{2} \frac{1 + \cos \pi y}{y}\ \de y\right) = \frac{I_0^2}{2c}\int_0^2 \frac{1 - \cos\pi y}{y}\ \de y$$
and finally with another easy change:
$$\mathcal{P} =   \frac{I_0^2}{2c}\int_0^{2\pi} \frac{1 - \cos\pi z}{z}\ \de z$$
\\
Making use of the relation on page \textcolor{Green}{22}, we can write
$$\mathcal{P} = \frac{I_0^2}{2c}\left[\gamma + \ln(2\pi) - \text{Ci}(2\pi)\right]$$
\newpage\noindent 

\section{Exercises}

\subsection*{Exercise 1.}

\textsf{Verify the following integral representation for the square of the Exponential Integral}:
$$[\text{Ei}(-z)]^2 = 2e^{-2z}\int_0^{+\infty} e^{-2zt}\frac{\ln(1+2t)}{1+t}\ \de t ~~~~~~~~~~~ |\arg(z)| \leq \frac{\pi}{2}$$
\\
\subsection*{Exercise 2.}

\textsf{Using L'Hôpital rule, show that}
$$\lim_{x\to +\infty} xe^{-x}\mathcal{E}_1(x) = 1$$
\\
\subsection*{Exercise 3.}

\textsf{Evaluating the following integral}
$$\int_{-\infty}^x \frac{e^t}{t^2(t-1)}\ \de t$$

\flushright
Solution: $\frac{e^x}{x} - 2\text{Ei}(x)  + e\text{Ei}(x-1)$, for $x<0$

\end{document}